\documentclass[12pt]{article}

\usepackage{amssymb,amsmath,latexsym}
\usepackage[dvips]{graphics}
\usepackage{mathrsfs}
\usepackage{float}

\newtheorem{prop}{Proposition}[section]

\newtheorem{lemma}{Lemma}[section]
\numberwithin{equation}{section}
\renewcommand{\baselinestretch}{1.4}
\def\ZZ{\mathbb{Z}}

\def\HH{\mathbb{H}}

\def\EE{\mathbb{E}}

\def\mmax{\mathrm{max}}
\def\mmin{\mathrm{min}}

\def\RR{\mathbb{R}}

\def\<{\langle}
\def\>{\rangle}
\def\pf{\noindent{\bf Proof.} }

\def\P{{\bf P}}

\def\qed{{\hfill $\Box$\medskip}}

\allowdisplaybreaks[4]

\input{epsf.sty}
\usepackage{epsfig}

\begin{document}
\title{\bf Phase Transition for the Contact Process in a Random Environment on $\ZZ^d\times\ZZ^+$}

\author{Qiang Yao\footnote{School of Statistics, East China Normal University and NYU--ECNU Institute of Mathematical Sciences at NYU Shanghai. E-mail: qyao@sfs.ecnu.edu.cn.}
        }

\maketitle{}

\begin{abstract}
We review the results in Chen \& Yao \cite{Chen-Yao2009}\cite{Chen-Yao2012} which concern the contact process in a static random environment on the half space $\ZZ^d\times\ZZ^+$ and make some addition to them. Furthermore, we explain why our methods cannot apply to the whole space case and compare our results with some related works.

 \end{abstract} \noindent{\bf 2000 MR subject
classification:} 60K35


\section{Introduction}

\subsection{Basic definitions of the contact process}

The basic contact process, which will be denoted by ``contact process'' in the following, was introduced in Harris \cite{Harris1974}. It is a model to describe  the spread of diseases. The process is defined as follows. Given a graph $G=(V,E)$, where $V$ denotes the vertex set of $G$, and $E$ denotes the edge set of $G$, the contact process $(\xi_{t}:~t\geq0)$ is a continuous--time Markov process, whose state space is $\{A:~A\subseteq V\}$. At each $t$, each vertex is either healthy or infected. Denote by $\xi_t$ the collection of infected vertices at time $t$. The transition rates are as follows:
\begin{equation}\label{e:transitionratesfix}
\left\{\begin{array}{ll} \xi_{t}\rightarrow\xi_{t}\setminus\{x\} \text{ for } x\in\xi_{t} \text{ at rate } 1,\\
\xi_{t}\rightarrow\xi_{t}\cup\{x\} \text{ for } x\notin\xi_{t}
\text{ at rate } \lambda\cdot|\{y\in\xi_{t}:x\sim
y\}|,\end{array}\right.
\end{equation}
where $\lambda>0$ is a positive constant, $|\cdot|$ denotes the cardinality of a set, and ``$x\sim y$'' denotes that the vertices $x$ and $y$ are neighbors. The intuitive interpretation of the above transition rates is that an infected vertex becomes healthy at fixed rate $1$, while a healthy vertex becomes infected at rate proportional to the number of its infected neighbors. The proportional coefficient $\lambda$ is the parameter of the contact process. Readers can refer to the standard references Liggett \cite{Liggett1985} and Durrett \cite{Durrett1988} for how the above rates determine a Markov process in a rigorous way. Often we use the notation $(\xi_t^A:~t\geq0)$ to denote the contact process with initial state $A$, that is, at time $0$ all vertices in $A$ are infected, while all vertices outside of $A$ are healthy. There is another viewpoint for the contact process which treats infected vertex as ``$1$'' while treats healthy vertex as ``$0$''. Under this viewpoint, the contact process is a Markov process with state space $\{0,1\}^V$. Therefore, the contact process is a special example of ``spin system''~(see Liggett \cite{Liggett1985} for rigorous definition).

The main problem in studying the contact process is its asymptotic behavior. For the process $\xi_t^O$ with a single infected vertex $O\in V$ at time $0$, we say
that the process survives if $\P(\xi^O_t\neq \emptyset~\text{for any}~t\geq0)>0$, otherwise we say that the process dies out. Furthermore, we say that the process survives strongly if $\P(\forall T\geq0,~\exists t>T,~\text{such that}~O\in\xi_t^O)>0$. And we say that the process survives weakly if it survives but not survives strongly. By the monotonicity~(or attractiveness) of the contact process~(which implies that the process is inclined to survive with larger infection parameter $\lambda$), we can define two critical values as follows:
$$\left\{\begin{array}{ll}\lambda_1:=\inf\{\lambda:~\xi_t^O~\text{survives}\},\\ \lambda_2:=\inf\{\lambda:~\xi_t^O~\text{survives strongly}\}. \end{array}\right.$$ If $G$ is a connected graph, then the value of $\lambda_1$ and $\lambda_2$ do not depend on the choice of the vertex $O$. Since strong survival implies survival, it can be easily seen that $\lambda_1\leq\lambda_2$.

\subsection{Known results for the contact process on $\ZZ^d$}

The contact process was firstly studied on the straight line $\ZZ^1$. Liggett \cite{Liggett1985} and Durrett \cite{Durrett1988} contain the main results for the one--dimensional case. The seminal work of Bezuidenhout \& Grimmett \cite{Bezuidenhout-Grimmett1990} used different geometric constructions to get the results for the high dimensional case, including:
\begin{itemize}
\item[\textbf{(a)}]~$\lambda_1=\lambda_2$~(denote by $\lambda_c$ the common value);
\item[\textbf{(b)}]~the process with parameter $\lambda_c$ dies out;
\item[\textbf{(c)}]~the complete convergence theorem holds for all $\lambda>0$, that is, for any $A\subseteq\ZZ^d$, $$\xi_t^A\Rightarrow
\overline{\nu}\cdot\P(\xi^A_t\neq \emptyset~\text{for any}~t\geq0)+\delta_{\emptyset}\cdot\P(\xi^A_t=\emptyset~\text{for some}~t\geq0)$$ as $t$ tends to infinity, where $\overline{\nu}$ denotes the upper invariant measure~(that is, the weak limit of the distribution of $\xi_t^{\ZZ^d}$ as $t\rightarrow\infty$), $\delta_\emptyset$ denotes the measure putting mass one on the empty set, and ``$\Rightarrow$'' stands for weak convergence;
\item[\textbf{(d)}]~the shape theorem holds, that is, there exists a convex subset $U\subseteq\RR^d$, such that for any $\varepsilon>0$, $$(1-\varepsilon)U\subseteq\frac{1}{t}H_t^0\subseteq(1+\varepsilon)U~~\text{eventually}$$ almost surely on the event that $\xi_t^0\neq\emptyset$ for any $t\geq0$, where $0$ denotes the origin of $\ZZ^d$, and $H_t^0=\displaystyle\bigcup\limits_{0\leq s\leq t}\xi_s^0$ denotes the set of vertices that have ever been infected before time $t$.\\
    \end{itemize}

\noindent\textbf{Remark.}~(1)~In the following, we use \textbf{(a)}, \textbf{(b)}, \textbf{(c)} and \textbf{(d)} to denote the above four results for short.

(2)~In Bezuidenhout and Grimmett \cite{Bezuidenhout-Grimmett1990}, the notations of the two critical values $\lambda_1$ and $\lambda_2$ have not been mentioned~(they first appeared in Pemantle \cite{Pemantle1992}). But Theorem 3 in \cite{Bezuidenhout-Grimmett1990} implies this conclusion.

(3)~Bezuidenhout and Grimmett \cite{Bezuidenhout-Grimmett1990} contains the proof of \textbf{(a)} and \textbf{(b)}. They didn't give the formal proof of \textbf{(c)} and \textbf{(d)}. The detailed proof of these results are provided in Liggett \cite{Liggett1999}.

\subsection{Contact processes in random environments}

Liggett \cite{Liggett1991} gives a general setting for the contact process in random environment. That is, the transition rates in (\ref{e:transitionratesfix}) are modified by
\begin{equation}\label{e:transitionratesrandom}
\left\{\begin{array}{ll} \xi_{t}\rightarrow\xi_{t}\setminus\{x\} \text{ for } x\in\xi_{t} \text{ at rate } \delta_x,\\
\xi_{t}\rightarrow\xi_{t}\cup\{x\} \text{ for } x\notin\xi_{t}
\text{ at rate } \sum\limits_{y\in\xi_{t},~y\sim
x}\lambda_{(y,x)},\end{array}\right.
\end{equation}
where $\{\delta_x:~x\in V\}$ and $\{\lambda_e:~e\in E\}$ are random variables chosen in a stationary ergodic manner. That means, the recovery rates and infection rates become random.\\

The contact process in random environment was first studied on $\ZZ^1$~(see Bramson, Durrett \& Schonmann \cite{Bramson-Durrett-Schonmann1991}, Liggett \cite{Liggett1991}\cite{Liggett1992}, Klein \cite{Klein1994}, Newman \& Volchan \cite{Newman-Volchan1996}, etc), focusing on the conditions for survival~(extinction). The high--dimensional case is more challenging. Chen \& Yao \cite{Chen-Yao2009}\cite{Chen-Yao2012} settled \textbf{(c)} in the half space case when $\delta_x\equiv1$ and $\lambda_e$'s are independent and identically distributed. Garet \& Marchand \cite{Garet-Marchand2012}\cite{Garet-Marchand2014} settled \textbf{(d)} when $\delta_x\equiv1$ and $\lambda_e$'s are stationary, ergodic and properly bounded.

All the above models belong to contact processes in static random environments, that is, the environment does not change as time goes. There are some models concerning contact processes in dynamic random environments; see, for example, Broman \cite{Broman2007}, Remenik \cite{Remenik2008}, Steif \& Warfheimer \cite{Steif-Warfheimer2008}, etc. The main difficulty in studying the processes in static random environments is that the process is not Markovian under the annealed~(or averaged) law.

\subsection{Organization of this article}

In Section 2, we will consider \textbf{(a)}, \textbf{(b)} and \textbf{(c)} in the half space case. \textbf{(c)} has been proved in Chen \& Yao \cite{Chen-Yao2012}, so we only state the proof sketch heuristically. \textbf{(a)} and \textbf{(b)} must be posed in a ``parameterized version'', and will be proved using the idea of Grimmett and Mastrand \cite{Grimmett-Mastrand1990}~(they considered the percolation model). A special case is an addition to Chen \& Yao \cite{Chen-Yao2009}, since we did not prove that the critical process dies out for the half space percolation cluster case there. In Section 3, we will compare our results with some related works and explain why our methods cannot apply to the whole space case.

\section{Contact Process in a Random Environment on $\ZZ^d\times\ZZ^+$}

The graph we are considering is $(\HH,\EE)$, where
$\HH=\ZZ^d\times\ZZ^+~(d\geq1)$, with $\ZZ=\{0,\pm 1,\pm
2,\cdots\}$ and $\ZZ^+=\{0,1,2,\cdots\}$; and $\EE=\{(x,y):~x,y\in\HH,~\|x-y\|=1\}$, with $\|\cdot\|$ denoting the
Euclidean norm. The graph is treated as unoriented; that is, $(x,y)$ and $(y,x)$ denote the same edge for all $x,y\in\HH$ satisfying $\|x-y\|=1$. The environment is defined via (\ref{e:transitionratesrandom}) with $\delta_x\equiv1$ and $\lambda_e$'s being i.i.d.$\sim\mu$, where $\mu([0,+\infty))=1$.

\subsection{The complete convergence theorem \textbf{(c)}}

The complete convergence theorem \textbf{(c)} was proved in Chen \& Yao \cite{Chen-Yao2009} for the half space percolation cluster case~(where $\mu$ follows the Bernoulli distribution), and was proved later in Chen \& Yao \cite{Chen-Yao2012} for the general half space case. In \cite{Chen-Yao2009} and \cite{Chen-Yao2012} we only proved the half plane case~(when $d=1$). The higher dimensional case~(when $d\geq2$) can be settled with no substantial difficulty.\\

\noindent\textbf{Step 1: The block conditions}

The ``building block'' of the proof procedure is the setup of the ``block conditions''~(denoted by \textbf{(BC)} in the following) for the survival of the process~(Proposition 3.2 in Chen \& Yao \cite{Chen-Yao2012}). Intuitively, suppose the process survives, then for any $\varepsilon>0$ sufficiently small, we can construct two kinds of boxes whose sizes depend on $\varepsilon$ but with almost fixed shape, such that with probability greater than $1-\varepsilon$, a horizontal seed~(i.e. an interval with all infected vertices) on the bottom of each box can give birth to another vertical seed with the same length on the right side of the box, with infection path being entirely contained in the interior of the box. The two kinds of boxes are called by ``S--box'' and ``L--box'' respectively, where ``S'' stands for ``short'' and ``L'' stands for ``long''. See Figure 1 for illustration. Note that when $d\geq2$, only one kind of boxes are needed by using the skew lines.

\begin{figure}[H]\label{picc:block}
 \center
 \includegraphics[width=10.0true cm, height=5true cm]{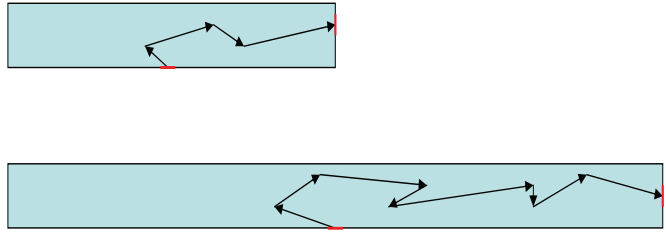}
 \caption{Construction of blocks}
 \end{figure}

The proof of Proposition 3.2 in Chen \& Yao \cite{Chen-Yao2012}~(which is the main contribution of that paper) was divided into three cases~(with three totally different proofs), one of which covers the proof of Lemma 3.4 in Chen \& Yao \cite{Chen-Yao2009} as a special case. Note that there is a similar disjunction in Garet \& Marchand \cite{Garet-Marchand2013}.\\

\noindent\textbf{Step 2: The dynamic renormalization construction}

Then we will use the S--boxes and L--boxes to construct a route, such that with probability greater than $1-\varepsilon$, a seed in a fixed square is joined through this route to some seeds in the other two fixed squares depending on $\varepsilon$ and having the same size (one above, the other on the right). See Figure 2 for illustration.

 \begin{figure}[H]\label{picc:uprightspread}
 \center
 \includegraphics[width=10true cm, height=7.5true cm]{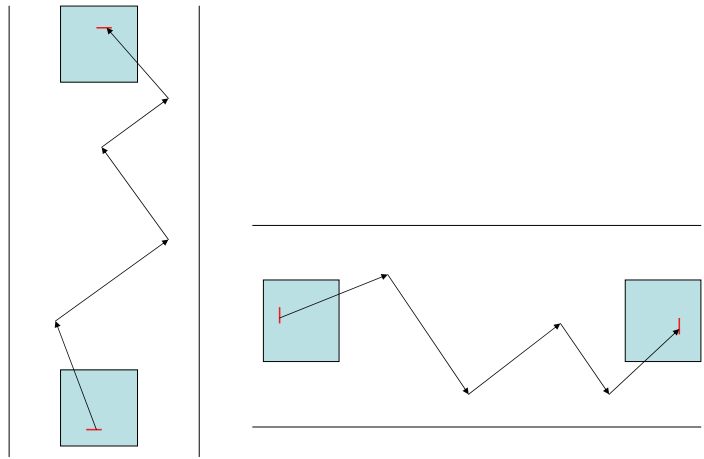}
 \caption{Producing new seeds in the north square and the east square}
 \end{figure}

Next, we fix $\varepsilon>0$ sufficiently small and iterate the above procedure several times in both directions~(to the ``east'' and to the ``north''), then treat the graph in a larger scale~(called the ``dynamic renormalization'' procedure). Figure 3 gives an illustration for the case that the ``large scale length'' of the ``large scale square'' is $3$. Furthermore, denote by $T(n,\varepsilon)$ the time span that the seed in the ``southwest'' generate the seed in the ``northeast'' in the ``large scale square'' with ``large scale length'' $n$.

 \begin{figure}[H]\label{picc:dynamicrenormalization}
 \center
 \includegraphics[width=10true cm, height=8true cm]{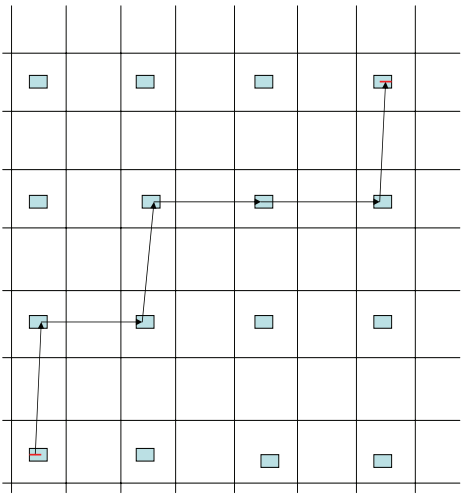}
 \caption{Dynamic renormalization}
 \end{figure}

The next proposition is Proposition 4.1 in Chen \& Yao \cite{Chen-Yao2012}, which is the main result for the dynamic renormalization. It tells us that as $n$ tends to infinity, $T(n,\varepsilon)$ follows the ``almost linear growth'' property asymptotically. The detailed proof can be found in Appendix 2 of Chen \& Yao \cite{Chen-Yao2009}.

\begin{prop}[Chen \& Yao \cite{Chen-Yao2012}]\label{p:dynamicrenormalization}Suppose that \textbf{(BC)} holds. Then there exists  $\overline{W}>0$, such that
$$\lim_{\varepsilon\rightarrow 0+}\liminf_{n\rightarrow\infty}\P\left( \frac{7\overline{W}}{6}n<T(n,\varepsilon)<\frac{11\overline{W}}{6}n\right)=1.$$\\
\end{prop}

\noindent\textbf{Step 3: Checking the equivalent conditions for (c)}

Having made the above preparations, we can prove \textbf{(c)} by checking the following two assertions~(Theorem 1.12 of Liggett \cite{Liggett1999}):
\begin{itemize}
\item[\textbf{(c1)}]~$\P\left(x\in\limsup\limits_{t\rightarrow\infty}\xi_t^A\right)=\P(\xi^A_t\neq \emptyset~\text{for any}~t\geq0)$ for all $x\in \HH$ and $A\subset\HH$.
\item[\textbf{(c2)}]~$\lim\limits_{M\rightarrow\infty}\liminf\limits_{t\rightarrow\infty}\P(\xi_t^{B_x(M)}\cap
B_x(M)\neq\emptyset)=1$ for all $x\in\HH$, where $B_x(M)$ is defined to be the ``ball'' centered at $x$ and with radius $M$~(but restricted on $\HH$).
\end{itemize}

\begin{figure}[H]\label{picc:CCTa}
\center
\includegraphics[width=12.5true cm, height=10true cm]{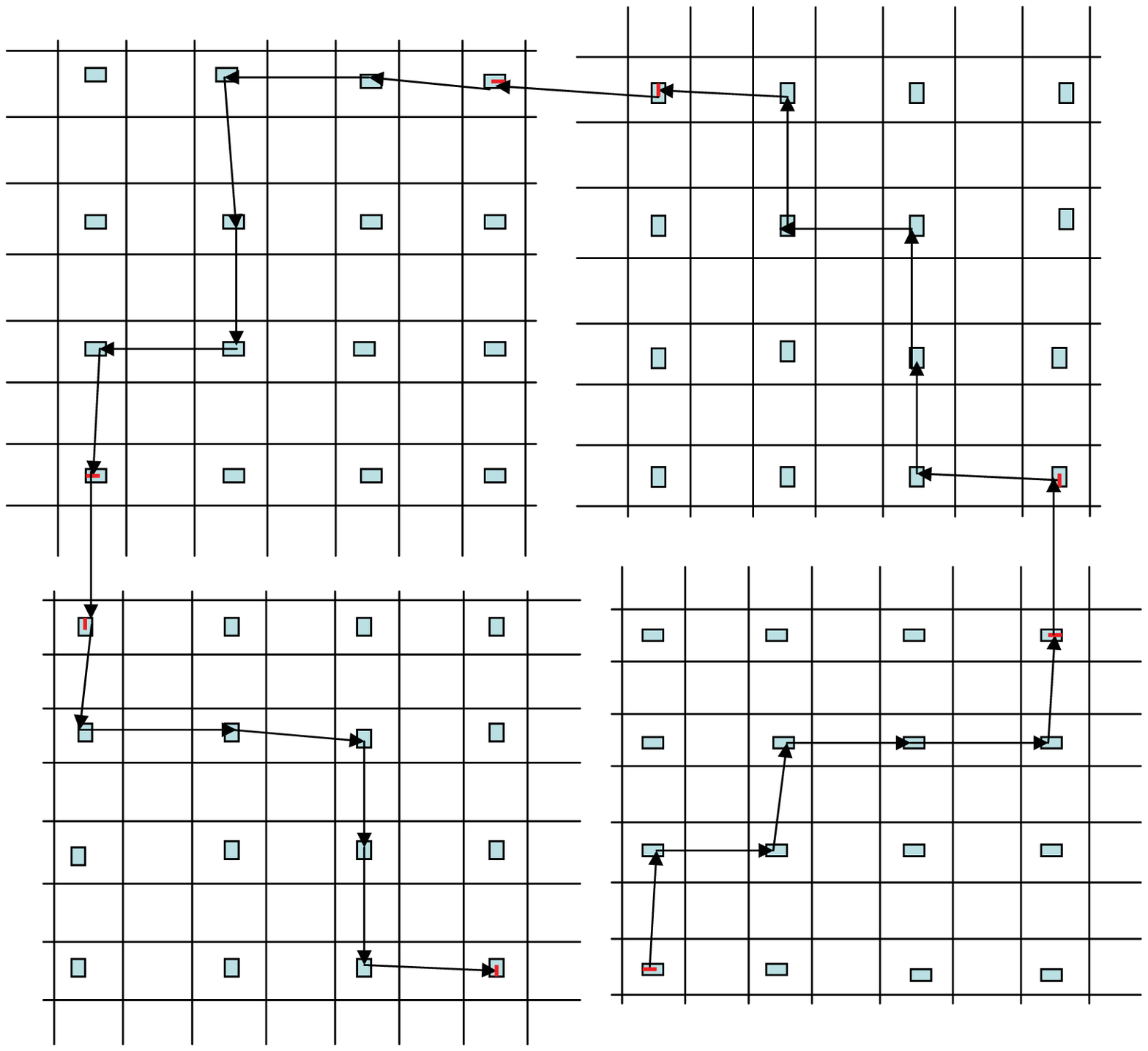}
\caption{Description of \textbf{(c1)}}
\end{figure}

The rigorous proof of \textbf{(c1)} and \textbf{(c2)} can be found in Subsections 5.1 and 5.2 of Chen \& Yao \cite{Chen-Yao2012}, respectively. The intuitive idea of \textbf{(c1)} is to iterate the construction posed in Proposition \ref{p:dynamicrenormalization} four times, as shown in Figure 4 intuitively. \textbf{(c2)} can be obtained by the above observation together with some extra tricks to prove the assertion that every remote site  cannot be infected in a short time.

\subsection{\textbf{(a)} and \textbf{(b)} in the parameterized version}

In order to consider \textbf{(a)} and \textbf{(b)}, we need to parameterize the model by changing the infection rate $\lambda_e$ by $\lambda\cdot\lambda_e$, where the $\lambda_e$'s are still i.i.d.~$\sim\mu$, and $\lambda>0$ is a free parameter. By monotonicity, we can still define $\lambda_1$ and $\lambda_2$, which are almost surely constants by translation invariance. Therefore, there is no difference between quenched law and annealed law when considering \textbf{(a)} and \textbf{(b)}. For simplicity, we use $\P$ for the measure, and use $\P_\lambda$ if we want to stress that the parameter is $\lambda$.\\

It is easy to see that \textbf{(c)} also holds for the parameterized version. Therefore, \textbf{(a)} holds trivially by \textbf{(c1)}. So we can denote by $\lambda_c$ the common value of $\lambda_1$ and $\lambda_2$.

\bigskip

\noindent\textbf{Remark.}~Whether $\lambda_c\in(0,+\infty)$ or not depends on the distribution $\mu$. When there exists $M\in(0,+\infty)$ such that $\mu([0,M])=1$, then $\lambda_c>0$. And when there exists $b\in(0,+\infty)$ such that $\mu([b,+\infty))=1$, then $\lambda_c<+\infty$. It will be interesting to consider the case when the support of $\mu$ is $(0,+\infty)$.\\

To prove \textbf{(b)}, that is, the critical process dies out in the parameterized version~(which includes the half space percolation cluster case in Chen \& Yao \cite{Chen-Yao2009}), we need the following lemma. The idea comes from Grimmett \& Mastrand \cite{Grimmett-Mastrand1990}, where they considered the percolation model.

\begin{lemma}\label{l:main}
If $\P_\lambda(\forall
t>0,~\xi_t^0\neq\emptyset)>0$, then there exists $\delta>0$, such that $$\P_{\lambda-\delta}(\forall t>0,~\xi_t^0\neq\emptyset)>0.$$
\end{lemma}

\pf Since $\P_\lambda(\forall
t>0,~\xi_t^0\neq\emptyset)>0$, the block conditions \textbf{(BC)} hold. Fix $\varepsilon>0$ sufficiently small as well as the variables~(including the sizes of the two kinds of boxes, the length of the seeds, and the time span) which guarantee \textbf{(BC)} to hold. Since all these variables have upper bounds depending only on the above $\varepsilon$, it follows from the continuity of the finite--time process in the parameter $\lambda$ that there exists $0<\delta<\lambda$, such that \textbf{(BC)} also hold with parameter $\lambda-\delta$ under the same $\varepsilon$ and the above variables. By Proposition \ref{p:dynamicrenormalization}, we have
$$\lim\limits_{\varepsilon\rightarrow0+}\liminf\limits_{n\rightarrow\infty}\P_{\lambda-\delta}(T(n,\varepsilon)<\infty)=1,$$ where $T(n,\varepsilon)$ is defined in Subsection 2.1. So we can choose $\varepsilon^\prime>0$ sufficiently small, such that $$\liminf\limits_{n\rightarrow\infty}\P_{\lambda-\delta}(T(n,\varepsilon^\prime)<\infty)>\frac{1}{2}.$$
Therefore, we have
$$\P_{\lambda-\delta}(T(n,\varepsilon^\prime)<\infty~~i.o.)\geq\liminf\limits_{n\rightarrow\infty}\P_{\lambda-\delta}(T(n,\varepsilon^\prime)<\infty)>\frac{1}{2}.$$
Furthermore, since $T(n,\varepsilon^\prime)<\infty~~i.o.$ implies that $\xi^{[-r,r]}$ can infect infinitely many sites, and therefore, the infection will persist forever. Here $r=r(\varepsilon^\prime)$ denotes half of the length of the initial seed in \textbf{(BC)} corresponding to the above $\varepsilon^\prime$. See Figure 4 for intuition. This implies
$$\P_{\lambda-\delta}(\forall t>0,~\xi_t^{[-r,r]}\neq\emptyset)>\frac{1}{2}.$$ Together with the trivial fact $\P_{\lambda-\delta}(\xi_1^0=[-r,r])>0$, we obtain
$$\P_{\lambda-\delta}(\forall t>0,~\xi_t^0\neq\emptyset)>0,$$
as desired.\qed\\

\noindent\textbf{Proof of \textbf{(b)}.}~Suppose $\P_{\lambda_c}(\forall t>0,~\xi_t^0\neq\emptyset)>0$. Then by Lemma \ref{l:main}, there exists $0<\delta<\lambda_c$, such that $\P_{\lambda_c-\delta}(\forall t>0,~\xi_t^0\neq\emptyset)>0$, contradicting with the definition of $\lambda_c$.\qed

\section{Concluding remarks and discussions}

The main idea of the proof procedure of \textbf{(c)}~(then \textbf{(a)} and \textbf{(b)} in the parameterized version) is enlightened by Bezuidenhout \& Grimmett \cite{Bezuidenhout-Grimmett1990}, that is, using the ``dynamic renormalization'' argument. The argument first appeared in Grimmett \& Mastrand \cite{Grimmett-Mastrand1990} and Barsky, Grimmett \& Newman \cite{Barsky-Grimmett-Newman1991}, where the authors considered the percolation model. But there are some big differences in our model. The ``block conditions'' in Bezuidenhout \& Grimmett \cite{Bezuidenhout-Grimmett1990} contain their Lemma 7~(which deals with ``space'' by using the fact that events depending on disjoint subgraphs are relatively independent) and Lemma 18~(which deals with ``time'' by using the Markov property of the process). However, In order to make good use of some symmetric properties, we need to consider the annealed law first~(Step 1 and Step 2 in the proof procedure), then go back to the quenched law to get the desired result~(Step 3 in the proof procedure). Under the annealed law, the fact that events depending on disjoint subgraphs are relatively independent still holds, but the Markov property does not hold any more. In consequence, if we consider the whole space case, we can only get a result similar to Lemma 7 in \cite{Bezuidenhout-Grimmett1990} and cannot get the result similar to Lemma 18 in \cite{Bezuidenhout-Grimmett1990}. And furthermore, we cannot get the desired result in the whole space case. On the other hand, the ``space block'' in $\ZZ^d\times\ZZ^+$ constructed by Proposition 3.2 in Chen \& Yao \cite{Chen-Yao2012} has similar function as the ``space--time block'' in $\ZZ^d\times\RR^+$ constructed by Lemmas 7 and 18 in Bezuidenhout \& Grimmett \cite{Bezuidenhout-Grimmett1990}. That is why we can get the results in the half space case.\\

We believe that the results are true for the whole space case. There may be some possible ways to prove the whole space case. The first possible idea is to prove directly. It is not easy. Even in the percolation cluster case, it is of the same difficulty as the long--existing problem that whether there is percolation at the critical point in the whole space case, which is clear in the half space case. The second possible idea is to prove that the critical value in the whole space case is the same as it in the half space case, which is clear for the percolation case as well as the contact process case, but it is not known for the contact process on the percolation cluster case. We will think about it in future research.\\

We mention at the end of this article that Garet \& Marchand \cite{Garet-Marchand2012}\cite{Garet-Marchand2014} deal with the shape theorem \textbf{(d)} in the whole space case under the assumption that $\delta_x\equiv1$ and $\lambda_e$'s are stationary, ergodic, and take value in $[\lambda_\mmin,\lambda_\mmax]$, where $\lambda_\mmin>\lambda_c(\ZZ^d)$, and $\lambda_\mmax<+\infty$. Getting rid of the boundedness assumption of $\lambda_e$'s may be a similar challenge as our model.\\

\bigskip

\noindent\textbf{Acknowledgment.}~We would like to thank Matthias Birkner, Rongfeng Sun, Jan Swart, as well as an anonymous referee for their useful suggestions on an earlier version of the manuscript. Partial supports from the Institute for Mathematical Sciences at the National University
of Singapore and a grant from the National Science Foundation of China~(No.11671145) are gratefully
acknowledged.

\end{document}